\documentclass{article}

\usepackage{graphicx}
\usepackage{amsmath}
\usepackage{amsthm}
\usepackage{amsfonts}
\usepackage{amssymb}
\usepackage{url}

\theoremstyle{definition}

\newcommand{\R}{\mathbb{R}}

\usepackage{caption}
\usepackage{subcaption}

\begin{document}

\title{A piecewise-linear isometrically immersed flat Klein bottle in Euclidean 3-space}
\author{Stepan Paul}

\maketitle

\begin{abstract}
 We present numerical polyhedron data for the image of a piecewise-linear map from a zero-curvature Klein bottle into Euclidean 3-space such that every point in the domain has a neighborhood which is isometrically embedded. To the author's knowledge, this is the first explicit piecewise-smooth isometric immersion of a flat Klein bottle. Intuitively, the surface can be locally made from origami and but for the self-intersections has the global topology of a Klein bottle.
\end{abstract}

\begin{figure}
 \centering
 \includegraphics[width=5in]{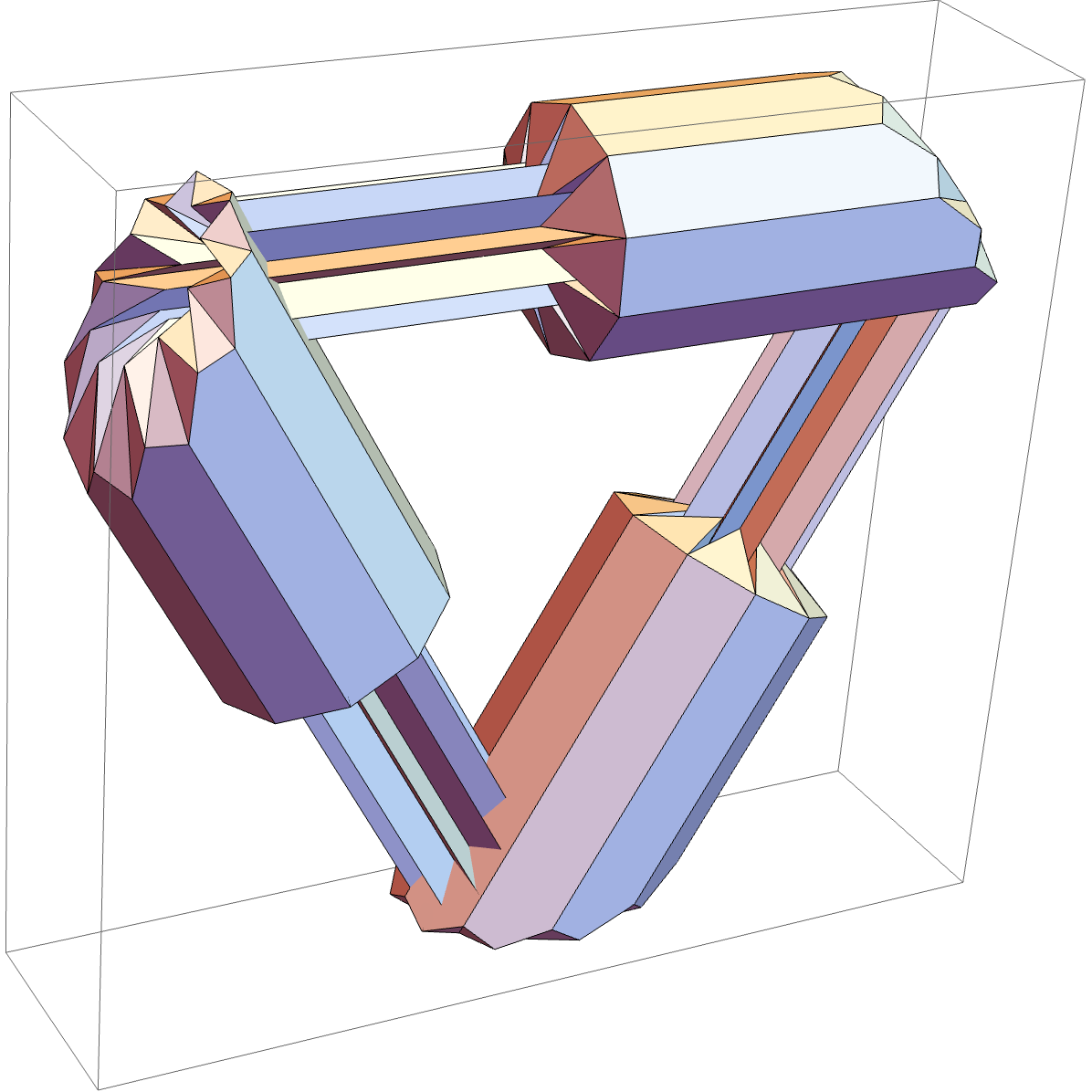}
 \caption{The image of a piecewise linear local isometric-embedding of a flat Klein bottle into $\mathbb{R}^3$.}
 \label{fig:fkb-10}
\end{figure}

The purpose of this manuscript is to announce a self-intersecting polyhedron, pictured in Figure \ref{fig:fkb-10}, that can be viewed as a piecewise-smooth isometric immersion of a zero-curvature Klein bottle.

 The polyhedron has 210 vertices, 510 edges, and 300 faces (210 triangles, 60 quadrilaterals, and 30 pentagons). The attached datasets include a combinatorial description of the polyhedral surface---given as a list of tuples in $\{1,2,\ldots,210\}$ defining the faces---and a list of points in $\R^3$ to which the vertices map. 
 
 One may check that the cell complex $\mathcal K$ resulting from the combinatorial data uses each edge exactly twice, showing that the complex defines a compact 2-manifold. Furthermore, the Euler characteristic is evidently zero, and one can check that the homology group has torsion, showing that $\mathcal K$ has the topology of a Klein bottle.
 
 One may also numerically check that each of the defined faces is planar, that each vertex has zero angle defect (i.e. the face angles at each vertex add to $2\pi$), and each vertex has an embedded vertex figure.
 
 The cell complex may be refined to a triangulation, so we may define a piecewise linear (with respect to barycentric coordinates on each triangle) map $\phi:\mathcal K\rightarrow\R^3$ sending vertices to their prescribed locations. Since a neighborhood of each vertex is embedded, $\phi$ is a local embedding, and so we can pull the polyhedral metric on the local image of $\phi$ back to $\mathcal K$. Then since each vertex has zero angle defect, this metric on $\mathcal K$ is flat in the sense of being locally isometric to a piece of the Euclidean plane.
 
  To the author's knowledge, this example furnishes the first explicit piecewise smooth isometric embedding of a flat Klein bottle into Euclidean 3-space. The existence of a piecewise linear map with these properties is implied by a result of Burago and Zalgaller \cite{bz1995}, although their proof is nonconstructive. This work builds upon \cite{me2021}---which exhibits a piecewise smooth, arc length preserving map from a flat Klein into $\R^3$ which fails to be locally injective at finitely many points---and runs parallel to constructions of piecewise linear embeddings of flat tori found, for example, in \cite{zalgaller2000,quintanar2020}.

\bibliographystyle{plain}
\bibliography{klein-bib}

\end{document}